# EV Dispatch Control for Supplementary Frequency Regulation Considering the Expectation of EV Owners

Hui Liu, *Member, IEEE*, Junjian Qi, *Member, IEEE*, Jianhui Wang, *Senior Member, IEEE*,
Peijie Li, *Member, IEEE*, Canbing Li, *Senior Member*, and Hua Wei

*Abstract*—Electric Vehicles (EVs) are promising to provide frequency regulation services due to their fast regulating characteristics. However, when EVs participate in Supplementary Frequency Regulation (SFR), it is challenging to simultaneously achieve the dispatch of the control center and the expected State of Charge (SOC) levels of EV batteries. To solve this problem, in this paper we propose a Vehicle-to-Grid (V2G) control strategy, in which an uncertain dispatch is implemented in the control center without detailed EV charging/discharging information. The regulation from the control center is achieved by allocating the regulation task within the frequency regulation capacity (FRC) of EVs. The expected SOC levels of EV batteries are guaranteed by a real-time correction of their scheduled V2G power in EV charging stations. Simulations on an interconnected two-area power system validate the effectiveness of the proposed V2G control in achieving both the frequency regulation and the expected SOC levels of EVs.

### NOMENCLATURE

| | |
|---|---|
| $B$ | Frequency bias factor |
| $E_i^{\text{rated}}$ | Rated capacity of the $i$th EV battery |
| $N_j$ | Number of EVs in the $j$th EV charging station |
| $p$ | Number of EV charging stations |
| $P_{i,k+1}^{\text{down}}$ | Regulation-down FRC of individual EV at time $k+1$ |
| $P_{i,k+1}^{\text{up}}$ | Regulation-up FRC of individual EV at time $k+1$ |
| $P_{i,k}$ | V2G power at the battery side of individual EV at time $k$ |
| $P_{\max}$ | Maximum V2G Power at the battery side of EVs |
| $P_i^{\text{sche}}(t)$ | Continuous scheduled V2G power at the battery side of individual EV for achieving the expected SOC |
| $P_{i,k+1}^{\text{sche}}$ | Scheduled V2G power at the battery side of individual EV for achieving the expected SOC at time $k+1$ |
| $P_i^{\text{const}}$ | Constant scheduled V2G power at the battery side of individual EV for achieving the expected SOC |
| $P_i^{\text{regu}}(t)$ | Continuous V2G regulation dispatch at the battery side of individual EV for performing SFR |
| $P_{i,k+1}^{\text{regu}}$ | V2G regulation dispatch at the battery side of individual EV for performing SFR at time $k+1$ |
| $R$ | A random ratio in the uncertain dispatch |
| $S_{j,k+1}^{\text{down}}$ | Regulation-down FRC of the $j$th EV charging station at time $k+1$ |
| $S_{k+1}^{\text{down}}$ | Total regulation-down FRC of all EVs at time $k+1$ |
| $S_{j,k+1}^{\text{up}}$ | Regulation-up FRC of the $j$th EV charging station at time $k+1$ |
| $S_{k+1}^{\text{up}}$ | Total regulation-up FRC of all EVs at time $k+1$ |
| $S_{k+1}^{\text{contr}}$ | Regulation task of EVs from the control center at time $k+1$ |
| $S_{k+1}^{\text{gener}}$ | Regulation task undertaken by conventional generating units at time $k+1$ |
| $S_{j,k+1}^{\text{station}}$ | Regulation task of the $j$th EV charging station at |

This work was supported in part by the National Natural Science Foundation of China (Grant No. 51577085) and in part, by the State Key Development Program for Basic Research of China (Grant No. 2013CB228205). J. Wang's work is supported by the U.S. Department of Energy (DOE)'s Office of Electricity Delivery and Energy Reliability.
H. Liu and H. Wei is with the Guangxi Key Laboratory of Power System Optimization and Energy Technology, College of Electrical Engineering, Guangxi University, Nanning, China (e-mail: hughlh@126.com; weihua@gxu.edu.cn).
J. Qi and J. Wang are with the Energy Systems Division, Argonne National Laboratory, Argonne, IL 60439 USA (e-mails: jqi@anl.gov; jianhui.wang@anl.gov).
P. Li is with the Guangxi Key Laboratory of Power System Optimization and Energy Technology, College of Electrical Engineering, Guangxi University, Nanning, China, and is visiting the Energy Systems Division, Argonne National Laboratory, Argonne, IL 60439 USA (e-mail: beyondpeijie@163.com).
C. Li is with the College of Electrical and Information Engineering, Hunan University, Hunan, China (e-mail: licanbing@qq.com).

| | |
|---|---|
| | time $k+1$ |
| $SOC_i^{exp}$ | Expected SOC of EV battery at plug-out time |
| $SOC_i^{initi}$ | Initial SOC of individual EV at plug-in time |
| $SOC_{i,k}$ | SOC of individual EV battery at time |
| $\Delta t^{regu}$ | Sampling time for frequency regulation |
| $\Delta t^{corr}$ | Sampling time for the scheduled V2G power |
| $t_i^{initi}$ | Plug-in time of an individual EV |
| $t_i^{out}$ | Plug-out time of an individual EV |
| $\eta^{ch}$ | Charging efficiency of EVs |
| $\eta^{disch}$ | Discharging efficiency of EVs |

## I. INTRODUCTION

ACCORDING to the Climate Change Synthesis Report, the share of low-carbon electricity supply will exceed 80% by 2050 and Renewable Energy Sources (RES) will be massively integrated into power grids [1]. Due to the random and intermittent nature of RES [2], [3] it is expected that the frequency excursion will become non-negligible. In order to maintain the system frequency in a tight band around its nominal value, the conventional generating units can be required to provide more frequency regulation capability. However, this will increase their regulation reserves, which is not cost-effective and may even incur undesirable high operational stresses.

As Electric Vehicles (EVs) can mitigate urban heat island effect to benefit local and global climates, the number of EVs continues to increase [4], [5].Therefore, EVs are promising to provide ancillary services to the power grid [6], [7], particular in frequency regulation, which can help reduce the frequency deviation and also the regulation reserve of the conventional generating units [8]. Recently, many Vehicle-to-Grid (V2G) control strategies have been proposed for both primary frequency control (PFC) and supplementary frequency regulation (SFR) [9]-[27].

The PFC of EVs is a decentralized V2G control that directly responds to the system frequency deviation. In [13] and [14], a constant droop control is developed to improve the frequency recovery performance. In [15], an aggregated model of EVs is proposed to evaluate the dynamic response in PFC. In addition, adaptive droop control methods have also been proposed [16]-[18]. By changing the droop coefficient according to the frequency deviation and the EV battery energy, the output power of EVs can be adjusted in real time to ensure both the regulation and the EV charging [18].

By contrast, the SFR of EVs is a dispatch-based V2G control for which the charging/discharging power of EVs is determined based on the regulation signal from the control center. Because the power provided by each EV is only up to 20 kW [19]-[21] while the power required by frequency regulation is in the order of megawatts [22], the aggregator has to be considered by the control center to perform SFR [23]. In addition, since EVs have their own transport usages, the expected State of Charge (SOC) levels of EV batteries should also be considered.

In [24] and [25], a V2G control is proposed to evaluate the performance of EVs in SFR. In [26], the controllable capacity of loads is used in V2G control to dispatch the regulation signal from the control center to EVs. However, the expected SOCs of EV owners are not thoroughly discussed [27]. In [28], a regulation dispatch is developed to perform SFR based on the day-ahead scheduling profile. However, it only addresses the frequency regulation but not the expected SOC of EV owners.

In [27], another V2G control for SFR is proposed to consider both the frequency regulation and the expected SOC levels of EVs. However, in order to coordinate the regulation-up and regulation-down dispatch, the expected charging/discharging power of the EV batteries has to be uploaded to the control center in real time. For this control strategy, the required information communication will inevitably increase the cost of V2G operation. In addition, in real power systems it is not practical for the system operators to consider the requirements of loads such as the charging demands of EVs.

In this paper, we propose a closed-loop V2G control to achieve both the frequency regulation and the EV charging demands. In particular, we develop a hierarchical control structure that consists of the control center, the EV aggregators, the EV charging stations, and individual EVs. Based on this control structure, an uncertain dispatch is implemented in the control center to achieve frequency regulation based on the frequency regulation capacity (FRC) of EVs and a real-time correction of the scheduled V2G power is performed in the charging stations to ensure the expected SOC levels of EVs.

The remainder of this paper is organized as follows. In Section III a hierarchical V2G framework is developed for SFR with EVs' participation. In Section III a dispatch control strategy of EVs is proposed. In Section IV simulation results are presented to validate the effectiveness of the proposed V2G control. Finally, conclusions are drawn in Section V.

## II. HIERARCHICAL FRAMEWORK OF ELECTRIC VEHICLES PARTICIPATING IN SFR

### A. SFR in Conventional Generation Control Systems

In power grid operation, the system frequency must be managed in a tight tolerance bound in order to maintain the supply-demand balance. For instance, in China the frequency deviation is usually kept within ±0.02 Hz. This is achieved by frequency regulation services such as the SFR (a centralized control performed in the control center).

The framework of SFR in the conventional generation system is illustrated in Fig. 1, in which the Area Control Error (ACE) reflects the supply-demand mismatch. When the Tie-line Bias Control (TBC) is considered as the operation mode of interconnected power grids, ACE is calculated based on the frequency and the tie-line power deviations. The main objective of SFR is to suppress the ACE fluctuation and keep the system frequency within the tolerance bound by adjusting the outputs of the generating units.

In order to implement the control shown in Fig. 1, the information that represents the characteristics of the generating units, such as the ramp speed, power outputs, upper regulation limits, and lower regulation limits, must be sent to the control center to be used for calculating the Frequency Regulation Capacity (FRC) and performing the Load Frequency Control (LFC). The outputs of the generating units will be adjusted according to the dispatch from the control center based on LFC.

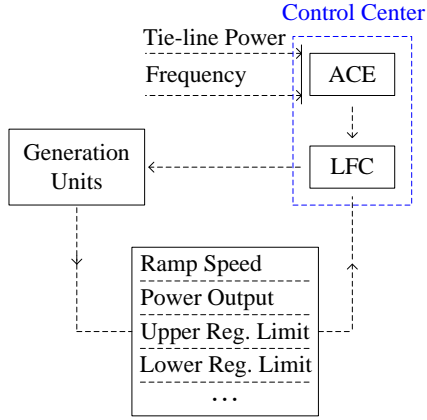

Fig. 1. The framework of SFR in the conventional generation control system

### B. Hierarchical V2G Framework of EVs Joining in SFR

As distributed energy storage systems, EVs can contribute to the SFR of interconnected power systems [24]-[28]. However, although the EV owners can choose to discharge their battery energy to the power grid to acquire revenues, the primary objective of EVs is for transport usages. Although [27] attempts to address the hierarchical V2G framework, the expected V2G power of the EV batteries has to be considered in the dispatch of the control center to ensure the expected battery SOC levels. In real power systems, it is not practical for the system operators to deal with the expected V2G power. Also, the dispatch in the control center is constrained by the expected V2G power to a great extent.

Here, a hierarchical framework with EVs participating in SFR is proposed, as illustrated in Fig. 2, which consists of four levels: 1) *Control Center level*; 2) *EV Aggregator level*; 3) *EV Charging Station level*; and 4) *Individual EV level*.

- In the Control Center level, the ACE is calculated based on the tie-line power and the frequency deviations and part of the ACE is randomly dispatched to the EV aggregators within the FRC of EVs.
- In the EV Aggregator level, the regulation dispatch from the Control Center is allocated to EV charging stations by the "V2G dispatch" block. The FRC uploaded by EV charging stations is summed by the "Total FRC" block.
- In the EV charging Station level, the regulation task is allocated to each EV by the "V2G Controller" block, and the charging power is regulated in real time to ensure the expected SOC levels of the EV batteries. The FRC of each EV is calculated by the "FRC Control" block. The "Information Management System" block communicates with individual EVs by the "Interface Circuit" block.
- In the Individual EV level, the "Interface Circuit" block controls the charging/discharging of an EV according to the command from EV charging station and uploads the EV information such as the real-time battery SOC, the plug-in duration, and the expected battery SOC.

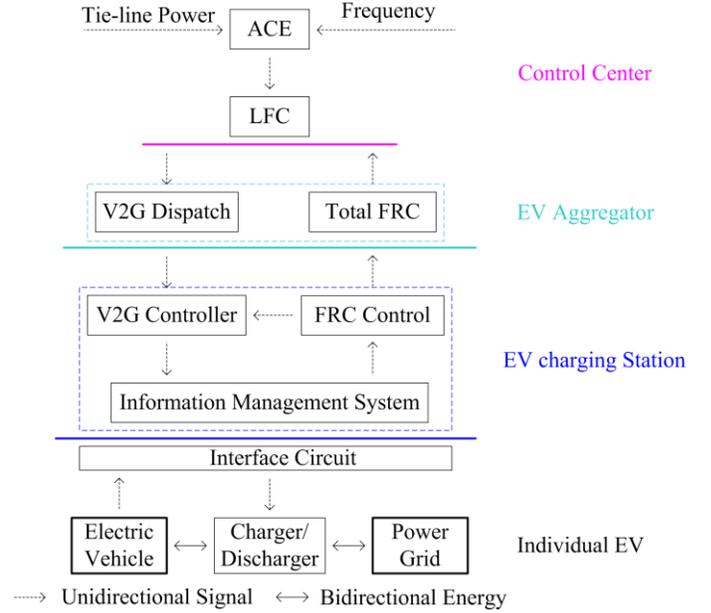

Fig. 2. The hierarchical framework for dispatch of EVs in SFR

### III. DISPATCH CONTROL OF EVs FOR PERFORMING SFR

Based on the hierarchical framework discussed above, a V2G control is proposed to perform SFR with EVs' participation.

#### A. Problem Description

As discussed in the introduction, the regulation dispatch from the control center and the expected battery energy levels of the EV owners are two concerns that have to be addressed.

**(1) Regulation dispatch from the control center**

Both ACE and FRC are crucial for the control center to implement the regulation dispatch to EVs, because ACE reflects the generation-load mismatch and FRC indicates the EVs' available capability for regulation. In conventional generation systems, the minimization of ACE is achieved only by the generating units. When EVs participate in SFR, they can undertake part of ACE for regulation. The dispatch of ACE to EVs should be performed within the FRC of EVs to make sure that the regulation task can be achieved. Also, the FRC of EVs has to be calculated in real time, because the V2G power that determines the FRC can vary with time due to the regulation and the scheduled charging.

Besides, the uncertainty of the regulation dispatch from the control center to EVs must be considered in the V2G control due to the uncertainty of the market, the randomness of EVs, the fluctuation of loads, etc.

**(2) Expected battery SOC levels**

The expected SOC levels of EV batteries are the most important concern of EV owners. Although the expected SOC levels of EVs are different from EV to EV due to different driving behaviors of EV owners, they can be categorized into three types: *increasing battery energy*, *decreasing battery energy*, and *maintaining battery energy*.

1) Increasing battery energy (called **TYPE I**): EV owners expect to charge their EVs to higher SOC levels for the next trip, especially when the battery energy of EVs is low. V2G control must ensure the expected battery levels before EV owners start the next travel.

2) Decreasing battery energy (called **TYPE II**): EV owners may want to get revenues by selling the redundant electricity to the power grid, especially when the battery energy is high enough for transport usages. The V2G control should be able to implement the discharging of EVs.

3) Maintaining battery energy (called **TYPE III**): Holding the battery energy is another choice, if EV owners do not want to charge their EVs or sell electricity to the power grid.

*B. Uncertain Dispatch in the Control Center*

In order to minimize the supply-demand mismatch as much as possible, the ACE as the objective of SFR is dispatched to the generating units in the automatic generation control system and also EVs. In particular, the ACE undertaken by EVs can be described as:

$$S_{k+1}^{\text{contr}} = f(\text{ACE}), \quad (1)$$

which represents the dispatch from the control center to EVs within their FRC. This function of ACE can depend on the market, the RES penetration, the reserve of generating units, the FRC of EVs, and many other factors. Therefore, the dispatch in (1) is uncertain. In order to explicitly describe this uncertainty, we rewrite (1) as

$$f(\text{ACE}) = R \cdot \begin{cases} -\min(|\text{ACE}|, S_k^{\text{up}}) & (\text{ACE} \leq 0) \\ \min(\text{ACE}, S_k^{\text{down}}) & (\text{ACE} > 0). \end{cases} \quad (2)$$

In (2), the dispatch depends on the ACE, the total regulation-up/regulation-down FRC ($S_k^{\text{up}}/S_k^{\text{down}}$) uploaded by the EV aggregator, and a random ratio $R$. Since the uncertain dispatch can depend on many factors, it is very difficult to accurately determine what distribution $R$ exactly follows. However, $R$ most probably follows the normal distribution mainly because 1) in real power systems the dispatches for regulation-down and regulation-up are usually close to each other [17], and 2) for a specific system there should be a most preferred proportion by which the EVs undertake frequency regulation in the control center and this proportion can be well represented by the mean value of a normal distribution. Therefore, we can further write $R$ as

$$R \sim N(\mu, \sigma^2) \quad (0 \leq R \leq 1), \quad (3)$$

where $\mu$ and $\sigma^2$ are, respectively, the mean and the variance of the normal distribution.

*C. Dispatch Control in EV Aggregators*

**(1) Total FRC**

In the conventional generation system, the reserve of generators, i.e. the FRC, is considered for regulation. Similarly, the FRC of EVs should also be calculated in order to properly perform the regulation task from the control center. For an EV aggregator, the total FRC can be calculated by

$$\begin{cases} S_{k+1}^{\text{up}} = \sum_{j=1}^{p} S_{j,k+1}^{\text{up}} \\ S_{k+1}^{\text{down}} = \sum_{j=1}^{p} S_{j,k+1}^{\text{down}} \end{cases}. \quad (4)$$

In (4), the total FRC of an EV aggregator is the summation of the FRC of all EV charging stations (i.e. $S_{j,k+1}^{\text{up}}/S_{j,k+1}^{\text{down}}$).

**(2) Regulation dispatch**

In an EV aggregator, the regulation task dispatched by the control center should be distributed to each EV charging station. Here the regulation task allocated to each EV charging station is proportional to its uploaded FRC:

$$S_{j,k+1}^{\text{station}} = \begin{cases} S_{k+1}^{\text{contr}} \cdot \dfrac{S_{j,k+1}^{\text{up}}}{S_{k+1}^{\text{up}}} & (S_{k+1}^{\text{contr}} \leq 0) \\ S_{k+1}^{\text{contr}} \cdot \dfrac{S_{j,k+1}^{\text{down}}}{S_{k+1}^{\text{down}}} & (S_{k+1}^{\text{contr}} > 0). \end{cases} \quad (5)$$

The rest of the regulation task is undertaken by the generating units as:

$$S_{k+1}^{\text{gener}} = \text{ACE} - S_{k+1}^{\text{contr}}. \quad (6)$$

*D. V2G Strategies in EV Charging Stations*

**(1) FRC calculation**

The FRC of an EV charging station at time $k+1$ can be calculated as

$$\begin{cases} S_{j,k+1}^{\text{up}} = \sum_{i=1}^{N_j} P_{i,k}^{\text{up}} \\ S_{j,k+1}^{\text{down}} = \sum_{i=1}^{N_j} P_{i,k}^{\text{down}} \end{cases} \quad (j=1,\cdots,p), \quad (7)$$

where

$$\begin{cases} P_{i,k}^{\text{up}} = P_{\max} + P_{i,k} \\ P_{i,k}^{\text{down}} = P_{\max} - P_{i,k}. \end{cases}$$

The FRC of each EV ($P_{i,k}^{\text{up}}/P_{i,k}^{\text{down}}$) is calculated based on the maximal V2G power $P_{\max}$ and the real-time V2G power $P_{i,k}$. $P_{\max}$ is a constant value depending on the charging/discharging devices while $P_{i,k}$ can change in real time due to the regulation and the scheduled charging. Therefore, the FRC of an EV charging station in (7) also varies with time.

The FRC for regulation-up and regulation-down provided by the individual EVs is illustrated in Fig. 3. When the V2G power is positive/negative, the regulation-up FRC will be greater/smaller than the regulation-down FRC. Therefore, the dispatches for regulation-up and regulation-down can be different due to the different FRC for regulation-up and regulation-down, which may result in deviation of the battery SOC from the expected. Thus the V2G power has to be adjusted to ensure the expected battery SOC levels of EV owners.

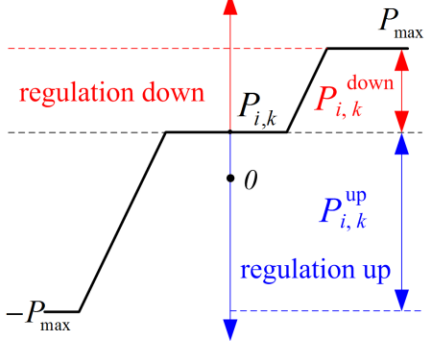

Fig. 3. FRC of individual EV at time $k$

**(2) V2G control**

The battery energy change of EVs comes from undertaking the regulation and performing the scheduled charging. The expected battery SOC of an EV can be described as:

$$\text{SOC}_i^{\text{exp}} = \text{SOC}_i^{\text{initi}} + \frac{\Delta E_i^{\text{regu}} + \Delta E_i^{\text{sche}}}{E_i^{\text{rated}}}, \quad (8)$$

where

$$\Delta E_i^{\text{regu}} = \int_{t_i^{\text{initi}}}^{t_i^{\text{out}}} P_i^{\text{regu}}(t)dt, \quad \Delta E_i^{\text{sche}} = \int_{t_i^{\text{initi}}}^{t_i^{\text{out}}} P_i^{\text{sche}}(t)dt.$$

In (8), the uncertain $\Delta E_i^{\text{regu}}$ resulted from regulation may lead to the deviation of battery SOC level from the expected. Therefore, a V2G control strategy must be developed to adjust $\Delta E_i^{\text{sche}}$ in order to achieve the expected SOC level of an EV.

If $\Delta E_i^{\text{regu}} = 0$, the constant scheduled V2G power of an EV can be calculated as:

$$P_i^{\text{const}} = \frac{(\text{SOC}_i^{\text{exp}} - \text{SOC}_i^{\text{initi}})E_i^{\text{rated}}}{t_i^{\text{out}} - t_i^{\text{initi}}}. \quad (9)$$

In (9), for TYPE I $P_i^{\text{const}} > 0$, for TYPE II $P_i^{\text{const}} < 0$, and for TYPE III $P_i^{\text{const}} = 0$.

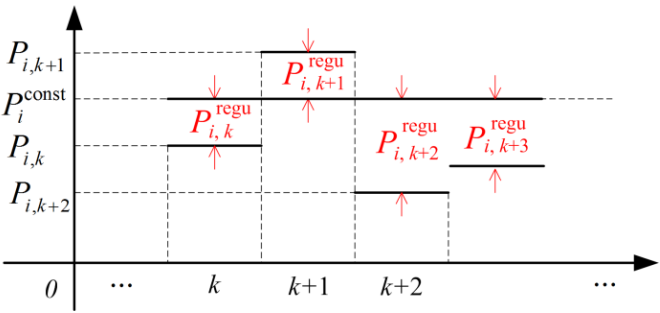

Fig. 4. V2G power with constant scheduled power and regulation dispatch.

The constant scheduled V2G power and the real-time dispatch are illustrated in Fig. 4. The V2G power of an EV will change around the constant scheduled V2G power due to regulation. If $\Delta E_i^{\text{regu}} = 0$, i.e., the energy change for regulation-up and regulation-down is equal, the expected battery SOC can be guaranteed by the constant scheduled V2G power. However, due to the uncertainty of the dispatch for the EV, $\Delta E_i^{\text{regu}}$ will probably not be zero. In that case, if the constant scheduled V2G power is still considered, the EV battery SOC will deviate from the expected level.

The dispatch decided by the control center is uncertain and uncontrollable by the EV. Therefore, in order to ensure the expected SOC level of an EV battery, the scheduled V2G power has to be adjusted in real time to compensate the change of the EV battery energy resulted from the uncertain dispatch. To achieve this, we propose a real-time scheduled V2G power based on (9) as:

$$P_{i,k}^{\text{sche}} = \frac{(\text{SOC}_i^{\text{exp}} - \text{SOC}_{i,q})E_i^{\text{rated}}}{t_i^{\text{out}} - t_q}. \quad (10)$$

Here $q = 1, 2, \cdots, t_q = t_i^{\text{initi}} + (q-1)\Delta t^{\text{corr}}$, and for each $q$, there is

$$k = (q-1)\frac{\Delta t^{\text{corr}}}{\Delta t^{\text{regu}}} + 1, \cdots, q\frac{\Delta t^{\text{corr}}}{\Delta t^{\text{regu}}}.$$

In (10), in order to make sure that $k$ is an integer, $\Delta t^{\text{corr}}$ is always chosen as an integer multiple of $\Delta t^{\text{regu}}$. The real-time scheduled V2G power of an EV is calculated based on the battery SOC at time $q$ and the remaining plug-in duration. It can respond to the change caused by frequency regulation in real time. The sample time interval of the real-time closed-loop control in (10) can be determined according to requirements.

For an individual EV, the regulation task at time $k+1$ can be allocated according to its FRC uploaded to the EV charging station at the previous time step as:

$$P_{j,k+1}^{\text{regu}} = \begin{cases} \dfrac{S_{j,k+1}^{\text{station}}}{\eta^{\text{disch}}} \cdot \dfrac{P_{j,k}^{\text{up}}}{\sum_{j=1}^{p} P_{j,k}^{\text{up}}} & (S_{j,k+1}^{\text{station}} \leq 0) \\ S_{j,k+1}^{\text{station}} \cdot \eta^{\text{ch}} \cdot \dfrac{P_{j,k}^{\text{down}}}{\sum_{j=1}^{p} P_{j,k}^{\text{down}}} & (S_{j,k+1}^{\text{station}} > 0). \end{cases} \quad (11)$$

In (11), the regulation task is assigned proportionally within the uploaded FRC to make sure that the regulation can be achieved by each EV. The effective charging or discharging efficiency may decrease if the battery is not charged or discharged at the prescribed rate. However, it is difficult to obtain a variable efficiency. Therefore, a constant efficiency is often used, as in [17] and [25].

With (10) and (11), the V2G power of an EV at time $k+1$ can be calculated by

$$P_{i,k+1} = P_{i,k+1}^{\text{sche}} + P_{j,k+1}^{\text{regu}}. \quad (12)$$

*E. Discussion on FRC and Scheduled V2G Power*

As in (12), if the regulation is not considered, the real-time

V2G power will increase when the scheduled V2G power increases. Therefore, the FRC of an EV decided by the real-time V2G power will be influenced by the scheduled V2G power. When the scheduled V2G power increases, the FRC for regulation up/down will increase/decrease. On the contrary, with the decrease of the scheduled V2G power, the FRC for regulation up decreases and that for regulation down increases.

If the scheduled V2G power is above or equal to the maximal V2G power, EVs will not be considered to participate in SFR. This is because in that case the EVs do not have the regulation down capability and the regulation up will result in further loss of battery energy. Therefore, charging will be the only option for an EV in that situation. Similarly, if the scheduled V2G power is below or equal to the negative maximal V2G power, EVs will not participate in SFR either.

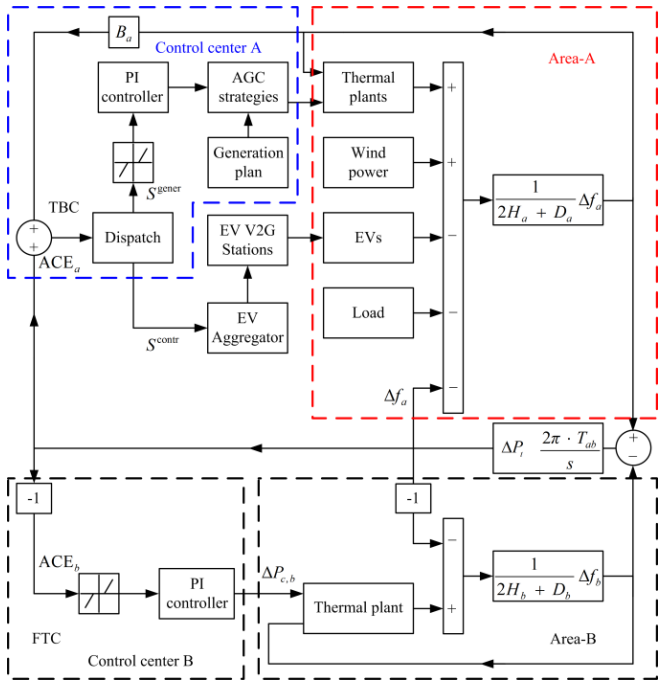

Fig. 5. A two-area interconnected power grid in China for SFR with EVs' participation

## IV. SIMULATION AND DISCUSSION

### A. Simulation System

As shown in Fig. 5, an interconnected two-area power system is modeled in MATLAB based on real power grid data in China [27] to simulate the SFR with EVs' participation.

The interconnected modes, the Tie-line Bias Control (TBC) and Flat Tie-line Control (FTC), are considered for area A and area B, respectively. We are most interested in area A since EVs and wind power are integrated here. In Fig. 6 we show the load profile, which comes from the historical data of a real power grid in China. In area A, 30 of 73 generating units take part in SFR, and the other generating units only follow generation curves. Besides, practical operation strategies are modeled into area A to imitate the SFR for conventional generating units. Basic parameters of this system are summarized in Table I.

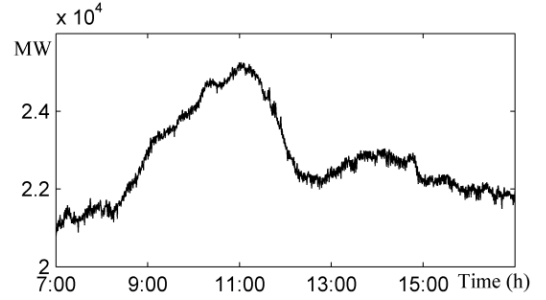

Fig. 6. Random load profile of a real power grid in China

The change of the EV battery energy is [18], [26].

$$\text{SOC}_{i,k} = \text{SOC}_i^{\text{initi}} + \frac{1}{E_i^{\text{rated}}} \Delta E_i, \quad (13)$$

where $\Delta E_i$ is the battery energy variation and satisfies

$$\Delta E_i = \int_0^k P_i(t)\,\mathrm{d}t. \quad (14)$$

TABLE I
PARAMETERS OF A TWO-AREA INTERCONNECTED POWER GRID

| Parameters | Area A | Area B |
| --- | --- | --- |
| Proportional and integral gains | 1, 0.01 | 1, 0.01 |
| Sampling time for frequency regulation (s) | 4 | 4 |
| Sampling time for the scheduled V2G power (h) | 1 | -- |
| Frequency bias factor (MW/0.1Hz) | 3400 | |
| Inertia constant (MW·s) | 16320 | 54720 |
| Load damping coefficient (MW/Hz) | 2040 | 6840 |
| Dead band of primary frequency control (Hz) | 0.033 | 0.033 |
| Communication delay (s) | 1 | 1 |
| Dead band of ACE (MW) | 20 | 20 |
| Time constant for wind power fluctuation (s) | 1 | 1 |

### B. Simulation Scenarios

#### (1) EV battery energy and driving behaviors

Different driving behaviors of EV owners can result in different residual battery energy levels of EVs at the parking lot. By assuming that the battery SOC levels of EVs follow the normal distribution shown in Table II [29], we simulate the battery SOC levels by Monte Carlo simulation method and build the battery models. In China, the EV users usually work from 8:00 to 17:00 [18], [27]. Thus, we assume that the plug-in duration of EVs is from 8:00 to 17:00.

TABLE II
RANDOM DISTRIBUTION OF DIFFERENT BATTERY SOC LEVELS OF EVS

|  | TYPE I | TYPE II | TYPE III |
|---|---|---|---|
| Initial SOC (pu) | SOC~$N$ (0.4,0.01) SOC $\in$ [0.3,0.5] | SOC~$N$ (0.7,0.01) SOC $\in$ [0.6,0.8] | SOC~$N$ (0.7,0.01) SOC $\in$ [0.6,0.8] |
| Expected SOC (pu) | SOC~$N$ (0.7,0.01) SOC $\in$ [0.6,0.8] | SOC~$N$ (0.4,0.01) SOC $\in$ [0.3,0.5] |  |

**(2) EV aggregators and charging stations**

An EV aggregator is assumed to be able to manage 100 EV charging stations, each of which serves 500 EVs. In each EV charging station, EVs have three types of behaviors, i.e., TYPE I, TYPE II, and TYPE III, as shown in Table III. In order to ensure the charging demands of EVs, an hourly recalculation of the scheduled V2G power is considered for each EV.

TABLE III
PARAMETERS OF AN EV CHARGING STATION

| EV Behaviors | Number of EVs | Total |
|---|---|---|
| TYPE I | 350 |  |
| TYPE II | 90 | 500 |
| TYPE III | 60 |  |

### C. V2G Control with Normally Distributed Dispatch

As discussed in Section III.B, the uncertain dispatch in the control center can be described by a normally distributed ratio $R \sim N(\mu,\sigma^2), 0 \leq R \leq 1$. Although for a specific system there should be a preferred proportion by which the EVs undertake frequency regulation in the control center, this proportion (actually the mean of the normal distribution $\mu$) can depend on many factors, such as the uncertainty of the market, the randomness of EVs, and the fluctuations of the RES. Determining such a proportion will require a lot of detailed information about the system which we do not have access to at this moment. Therefore, here we only choose $\mu = 0.5$ as an example for which the control center prefers to dispatch half of the ACE to EVs. In Section III.D we will further discuss how different mean values could influence the control performance.

The V2G control strategy with the constant scheduled V2G power in (9) is called "CS1" and that with the real-time scheduled V2G power in (10) is called "CS2". The case in which EVs do not participate in SFR is called "W/O V2G". In order to quantify the effectiveness of the V2G control, the Maximum values (Max), Minimum values (Min), and Root Mean Square values (RMS) of the ACE and the frequency deviation are calculated and listed in Tables IV and V.

TABLE IV
QUALITY OF THE ACE WITH NORMALLY DISTRIBUTED DISPATCH IN AREA A FOR $R \sim N$ (0.5, 0.01)

| Strategies | Max (MW) | Min (MW) | RMS (MW) |
|---|---|---|---|
| W/O V2G | 314.29 | -419.11 | 127.35 |
| CS1 | 298.86 | -341.48 | 88.13 |
| CS2 | 306.20 | -335.82 | 88.61 |

TABLE V
QUALITY OF FREQUENCY DEVIATION WITH NORMALLY DISTRIBUTED DISPATCH IN AREA A FOR $R \sim N$ (0.5, 0.01)

| Strategies | Max (MW) | Min (MW) | RMS (MW) |
|---|---|---|---|
| W/O V2G | 0.0676 | -0.0818 | 0.0252 |
| CS1 | 0.0599 | -0.0633 | 0.0175 |
| CS2 | 0.0550 | -0.0606 | 0.0175 |

With EVs participating in the frequency regulation, both CS1 and CS2 can better suppress the ACE and the frequency fluctuations when compared with the W/O V2G case. This is mainly because EVs have faster regulating characteristics than the conventional generating units.

Besides, we present the battery SOC of a randomly chosen EV in each type of EVs for illustration. The results for the other EVs are similar and thus are not presented. As shown in Figs. 7-9, CS2 can guarantee the expected battery SOC levels of EVs while CS1 cannot. This is because as shown in Figs. 10-12, CS1 uses the constant scheduled V2G power, while CS2 performs a real-time correction of the scheduled V2G power. Without real-time correction CS1 cannot compensate the battery energy change resulted from the uncertain regulation.

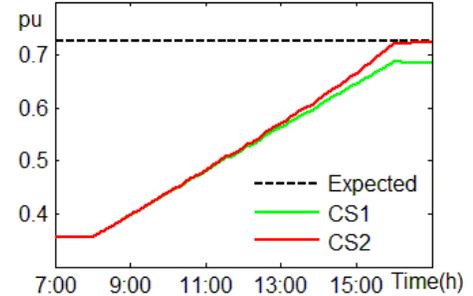
Fig. 7. Real-time battery SOC levels of an EV in TYPE I for $R \sim N$ (0.5, 0.01)

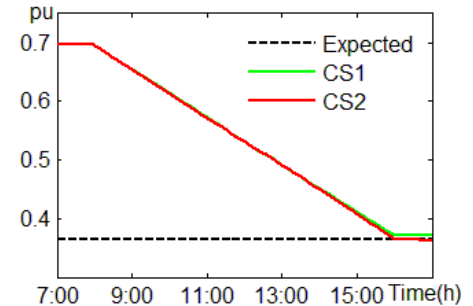
Fig. 8. Real-time battery SOC levels of an EV in TYPE II for $R \sim N$ (0.5, 0.01)

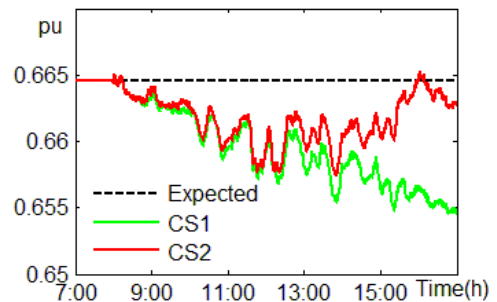
Fig. 9. Real-time battery SOC levels of an EV in TYPE III for $R \sim N$ (0.5, 0.01)

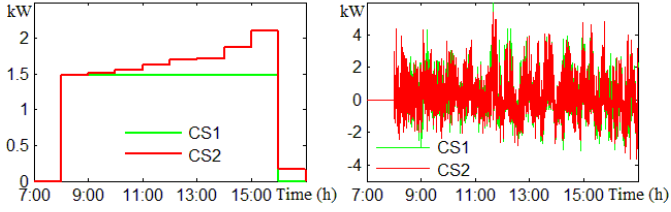

(a) Scheduled V2G power    (b) Regulation dispatch
Fig. 10. V2G power of an EV in TYPE I for CS1 and CS2

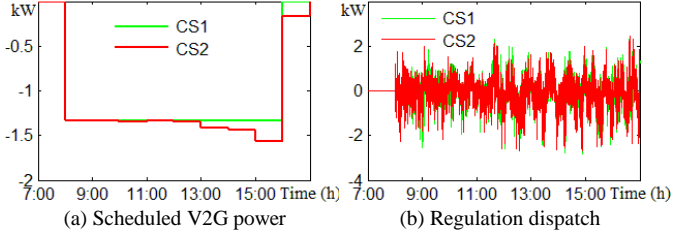

(a) Scheduled V2G power    (b) Regulation dispatch
Fig. 11. V2G power of an EV in TYPE II for CS1 and CS2

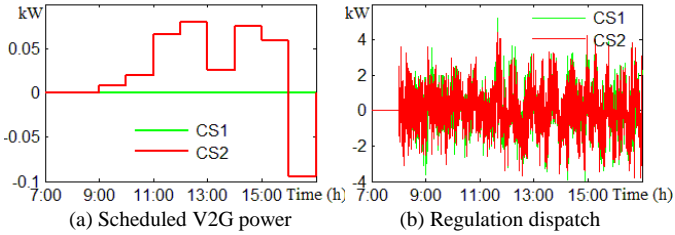

(a) Scheduled V2G power    (b) Regulation dispatch
Fig. 12. V2G power of an EV in TYPE III for CS1 and CS2

In order to more clearly show the effectiveness of CS2 in achieving the charging demands of EVs, we randomly choose 50 EVs with the same expected SOC levels in TYPE I and TYPE II. As shown in Fig. 13, the charging demands of all these EVs can be achieved by CS2.

In addition, it is seen in Fig. 13 and Figs. 7-8 that the charging/discharging curves have nearly constant slopes. This is not because of a relatively ideal load profile since there are significant fluctuations in the load profile shown in Fig. 6. This is actually due to the following two reasons. First, the EV power is very low and the charging/discharging in very short time periods will not result in a significant change of EV battery SOC levels. Second, there is no big gap between the regulation-up and the regulation-down.

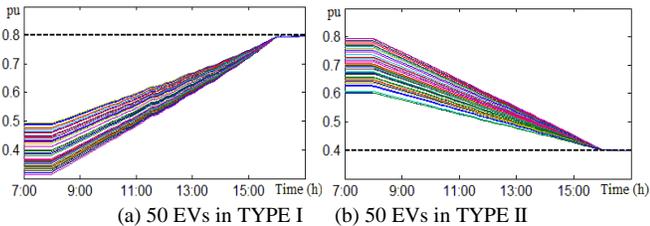

(a) 50 EVs in TYPE I    (b) 50 EVs in TYPE II
Fig. 13. Real-time battery SOC levels (The black dashed line represents the expected SOC level and the other curves represent the real time SOC levels)

### D. Control Performance under Different Mean Values of R

The mean value of the normal distribution for $R$ determines what proportion of frequency regulation that the EVs will undertake. Here we discuss how different mean values influence the control performance. In particular, we simulate another case in which $R \sim N(0.8, 0.01)$ and the corresponding Max, Min, and RMS of the ACE and frequency deviations are listed in Tables VI and VII, respectively. Comparing the results in Tables VI and VII and those in Tables IV and V, we can see that the dispatch with a greater mean value of $R$ has a better performance in suppressing the ACE and the frequency deviations. This is because under bigger mean value of $R$ more dispatch is distributed to EVs which have faster regulating characteristics than the conventional generating units.

TABLE VI
QUALITY OF THE ACE WITH NORMALLY DISTRIBUTED DISPATCH IN AREA A FOR $R \sim N$ (0.8, 0.01)

| Strategies | Max (MW) | Min (MW) | RMS (MW) |
| --- | --- | --- | --- |
| W/O V2G | 314.29 | -419.11 | 127.35 |
| CS1 | 311.94 | -417.42 | 84.19 |
| CS2 | 295.96 | -422.94 | 83.33 |

TABLE VII
QUALITY OF FREQUENCY DEVIATION WITH NORMALLY DISTRIBUTED DISPATCH IN AREA A FOR $R \sim N$ (0.8, 0.01)

| Strategies | Max (MW) | Min (MW) | RMS (MW) |
| --- | --- | --- | --- |
| W/O V2G | 0.0676 | -0.0818 | 0.0252 |
| CS1 | 0.0571 | -0.0772 | 0.0169 |
| CS2 | 0.0545 | -0.0748 | 0.0167 |

### E. Robustness to Different Distributions of R

As discussed in Section III-B, the ratio $R$ in the uncertain dispatch most probably follows a normal distribution. A uniform distributed $R \sim U[0,1], 0 \leq R \leq 1$ may not be a very realistic dispatch since the uncertain dispatch will probably not be so random. However, the uniform distribution can be used to test the robustness of the proposed V2G control to different uncertain dispatches in the control center. In Tables VIII and IX, we list the Max, Min, and RMS of the ACE and frequency deviations for $R \sim U[0,1], 0 \leq R \leq 1$. Comparing the results in Tables VIII and IX and those in Tables IV and V, we can see that when $R \sim U[0,1]$ the ACE and the frequency quality is not as good as that for $R \sim N(0.5, 0.01)$. This is mainly because $R \sim U[0,1]$ has much greater variance than $R \sim N(0.5, 0.01)$ although they have the same mean value. When the variance is greater, it means that there will be a higher probability to have an $R$ that is much smaller or much higher than the mean value. As is discussed in Section IV.D, greater $R$ is preferable since the EVs have faster regulating characteristics than conventional generators. Then greater variance leading to worse control performance seems to indicate that the impact from a smaller $R$ such as 0.3 that deteriorates the performance is higher than that from a greater $R$ such as 0.7 that improves the performance.

When we compare the results in Tables VIII and IX and those in Tables VI and VII, it is seen that the V2G control with $R \sim N(0.8, 0.01)$ has much better performance in suppressing the ACE and the frequency deviations than that with $R \sim U[0,1]$. This is because the former distribution has greater mean value and smaller variance.

As shown in Tables VIII and IX, when the same mean values and variances are considered for the two dispatches, the results for $R \sim U[0,1]$ and $R \sim N(0.5, 1/12)$ are very close to each other, indicating that the control performance of the proposed control mainly depends on the mean value and the variance of $R$ and the proposed approach is robust to the distributions of $R$.

TABLE VIII
QUALITY OF THE ACE WITH NORMALLY DISTRIBUTED DISPATCH IN AREA-A FOR $R \sim U[0, 1]$ AND $R \sim N(0.5, 1/12)$

| Random Ratio | Strategies | Max (MW) | Min (MW) | RMS (MW) |
|---|---|---|---|---|
|  | W/O V2G | 314.29 | -419.11 | 127.35 |
| $R \sim U[0, 1]$ | CS1 | 319.86 | -327.75 | 91.26 |
|  | CS2 | 309.02 | -330.38 | 91.18 |
| $R \sim N(0.5, 1/12)$ | CS1 | 308.23 | -317.31 | 91.24 |
|  | CS2 | 296.73 | -335.98 | 90.46 |

TABLE IX
QUALITY OF FREQUENCY DEVIATION WITH NORMALLY DISTRIBUTED DISPATCH IN AREA-A FOR $R \sim U[0, 1]$ AND $R \sim N(0.5, 1/12)$

| Random Ratio | Strategies | Max (Hz) | Min (Hz) | RMS (Hz) |
|---|---|---|---|---|
|  | W/O V2G | 0.0676 | -0.0818 | 0.0252 |
| $R \sim U[0, 1]$ | CS1 | 0.0617 | -0.0682 | 0.0188 |
|  | CS2 | 0.0611 | -0.0691 | 0.0187 |
| $R \sim N(0.5, 1/12)$ | CS1 | 0.0584 | -0.0602 | 0.0187 |
|  | CS2 | 0.0618 | -0.0630 | 0.0184 |

*F. Robustness to Initial Battery SOC Levels*

In the above discussions, the initial battery SOC levels of EVs are assumed to follow a normal distribution with a variance of 0.01, as shown in Table II. In real world, the initial battery SOC levels may have a bigger variance. Therefore, we also test a case in which the variance of the normal distribution is increased to 0.1 and the initial battery SOC is bounded in [0.1, 0.9] for all three types of EVs.

The RMS of the differences between the battery SOC levels and the expected during 16:00 to 17:00 is calculated. For both cases with different variances, the results for the EVs in each type are very close to each other. For demonstration, we randomly choose one EV in TYPE I, TYPE II, and TYPE III. As shown in Table X, when the variance of the initial battery SOC levels is significantly increased to 0.1, the RMS is still very small. Therefore, the proposed CS2 control is very robust to the initial battery SOC of EVs.

TABLE X
RMS OF THE DEVIATIONS BETWEEN THE BATTERY SOC AND THE EXPECTED

| Variances | RMS (pu) | | |
|---|---|---|---|
|  | TYPE I | TYPE II | TYPE III |
| 0.01 | 0.00058 | 0.00098 | 0.00076 |
| 0.1 | 0.00053 | 0.00097 | 0.00069 |

V. CONCLUSION

In this paper, the SFR with EVs' participation is achieved by a proposed V2G control that considers both the regulation from the control center and the expected battery SOC levels of the EV owners. Similar to the SFR of the conventional generation control system, a hierarchical V2G control framework is proposed for EVs to participate in SFR, which includes the control center, the EV aggregator, the EV charging station, and the individual EVs. An uncertain dispatch is considered in the control center without detailed EV charging/discharging information. Both the regulation and the expected charging of EVs can be ensured, because the regulation task is allocated within the FRC of EVs and the scheduled V2G power is corrected in real time.

The effectiveness of the proposed V2G control is validated by simulations on an interconnected two-area power grid model. The results show that the ACE and the frequency fluctuations of the interconnected power grids can be well suppressed when EVs participate in frequency regulation. At the same time, the expected SOC of EV batteries can also be guaranteed by the proposed V2G control.

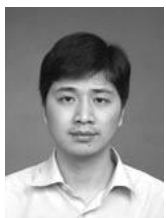
**Hui Liu** (M'12) received the M. S. degree in 2004 and the Ph.D. degree in 2007 from the College of Electrical Engineering, Guangxi University, Nanning, China, both in electrical engineering.

He was a Postdoctoral Fellow at Tsinghua University, Beijing, China, from 2011 to 2013 and was a staff at Jiangsu University, Zhenjiang, China, from 2007 to 2016. He visited the Energy Systems Division at Argonne National Laboratory, Argonne, IL, USA, as a visiting scholar from 2014 to 2015. He joined the College of Electrical Engineering at Guangxi University, Nanning, China, in 2016, where he is a professor. His research interests include power system control, electric vehicles, demand response, and power system optimization.

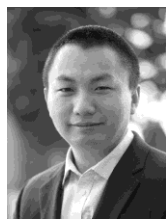
**Junjian Qi** (S'12–M'13) received the B.E. degree from Shandong University, Jinan, China, in 2008 and the Ph.D. degree from Tsinghua University, Beijing, China, in 2013, both in electrical engineering.

In February–August 2012 he was a Visiting Scholar at Iowa State University, Ames, IA, USA. During September 2013–January 2015 he was a Research Associate at Department of Electrical Engineering and Computer Science, University of Tennessee, Knoxville, TN, USA. Currently he is a Postdoctoral Appointee at the Energy Systems Division, Argonne National Laboratory, Argonne, IL, USA. His research interests include cascading blackouts, power system dynamics, state estimation, synchrophasors, electric vehicles, and cybersecurity.

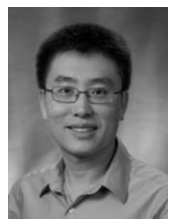
**Jianhui Wang** (M'07-SM'12) received the Ph.D. degree in electrical engineering from Illinois Institute of Technology, Chicago, IL, USA, in 2007.

Presently, he is the Section Lead for Advanced Power Grid Modeling at the Energy Systems Division at Argonne National Laboratory, Argonne, IL, USA. Dr. Wang is the secretary of the IEEE Power & Energy Society (PES) Power System Operations, Planning & Economics Committee. He is an associate editor of Journal of Energy Engineering and an editorial board member of Applied Energy. He is also an affiliate professor at Auburn University and an adjunct professor at University of Notre Dame. He has held visiting positions in Europe, Australia and Hong Kong including a VELUX Visiting Professorship at the Technical University of Denmark (DTU). Dr. Wang is the Editor-in-Chief of the IEEE Transactions on Smart Grid and an IEEE PES Distinguished Lecturer. He is also the recipient of the IEEE PES Power System Operation Committee Prize Paper Award in 2015.

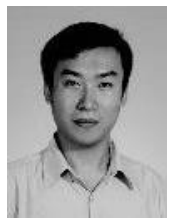
**Peijie Li** (M'15) received the B.E. and Ph.D. degrees in electrical engineering from Guangxi University, Nanning, China, in 2006 and 2012, respectively.

Since 2015, he has been a visiting scholar with the Argonne National Laboratory, Lemont, IL, USA. He is currently an Associate Professor at Guangxi University. His research interests include optimal power flow, small-signal stability, and security constrained economic dispatch and restoration.

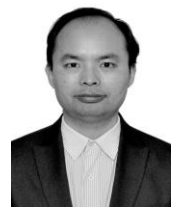
**Canbing Li** (M'06-SM'13) received the B.E. and Ph.D. degree from Tsinghua University, Beijing, China, in 2001 and 2006, respectively, both in electrical engineering.

Since 2014, he is currently a Professor at Hunan University, Changsha, China. His research interests include micro-grid, big data for energy systems, and power systems control.

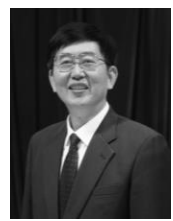
**Hua Wei** received the B.S. and M.S. degrees in power engineering from Guangxi University, China, in 1981 and 1987, respectively, and the Ph.D. degree in power engineering from Hiroshima University, Japan, in 2002.

From 1994 to 1997, he was a Visiting Professor with Hiroshima University, Japan. From 2004 to 2014, he was the Vice President of Guangxi University. He is currently a Professor of Guangxi University. He is also the Director of the Institute of Power System Optimization, Guangxi University. His research interests include power system operation and planning, particularly in the application of optimization theory and methods to power systems.